\documentclass[12pt]{article}
\usepackage{amsfonts}
\usepackage{amsmath,amssymb,amsthm}
\usepackage[usenames]{color}
\usepackage{mathrsfs}
\usepackage{syntonly}
\usepackage[top=1in,bottom=1.4in,left=1in,right=1in]{geometry}
\usepackage[all,cmtip]{xy}
\usepackage{bbm}
\usepackage{mathrsfs}
\usepackage{dsfont}
\usepackage{enumitem}
\usepackage{fouridx}
\usepackage[stable]{footmisc}
\usepackage{array}
\usepackage{caption}
\usepackage{latexsym}
\usepackage{extarrows}
\usepackage{amsfonts}
\usepackage{amsmath}
\usepackage{amsthm}
\usepackage{amssymb}
\usepackage{latexsym}
\usepackage{indentfirst}
\allowdisplaybreaks

\numberwithin{equation}{section}
\newtheorem{theorem}{Theorem}[section] %to number theorems, etc. within  sections
\newtheorem{proposition}[theorem]{Proposition}
\newtheorem{corollary}[theorem]{Corollary}
\newtheorem{lemma}[theorem]{Lemma}
\newtheorem{definition}[theorem]{Definition}
\newtheorem{remark}[theorem]{Remark}

\begin{document}

%%%  Abbreviations

\title{\Large\bf  Local Weyl modules and cyclicity of tensor products for  Yangians of $G_2$}
\author{{Yilan Tan}}

\date{} % 28 June 2005 - start.

\maketitle

\vspace{17 mm}
\begin{abstract}
Let $\mathfrak{g}$ be the exceptional complex simple Lie algebra of type $G_2$. We provide a concrete cyclicity condition for the tensor product of fundamental representations of the Yangian $Y(\mathfrak{g})$. Using this condition, we show that every local Weyl module is isomorphic to an ordered tensor product of fundamental representations of $Y(\mathfrak{g})$.
\end{abstract}
{\it Key words: Yangian; $G_2$; Local Weyl modules; Cyclicity condition}
\newpage
\section{Introduction}
There is a rich structure theory behind the finite-dimensional representations of Yangians $Y(\mathfrak{g})$ and quantum affine algebras $U_q(\hat{\mathfrak{g}})$ as the category of their finite-dimensional representations is not semi-simple, where $\mathfrak{g}$ is a complex simple Lie algebra of rank $l$.
Finite-dimensional irreducible representations of $U_q(\hat{\mathfrak{g}})$  are parameterized by $l$-tuples of polynomials $\operatorname{P}=\big(P_1(u),\ldots, P_l(u)\big)$, where $P_i(0)=1$, see \cite{ChPr6}. In \cite{ChPr1}, V. Chari and A. Pressley showed that in the class of all highest weight representations associated to $\operatorname{P}$, there is a unique (up to isomorphism) finite-dimensional highest weight representation $W(\operatorname{P})$ such that any other representation in this class is a quotient of $W(\operatorname{P})$. It has been established that $W(\operatorname{P})$ is isomorphic to an ordered tensor product of fundamental representations of $U_q(\hat{\mathfrak{g}})$, and a proof of this fact can be found in \cite{ChMo2}. In \cite{Ch3}, V. Chari provided a method to find a concrete cyclicity condition for an ordered tensor product of Kirillov-Reshetikhin modules using braided group action on the imaginary root vectors.

The finite-dimensional representation theory of $Y(\mathfrak{g})$ is an analogue of the one of $U_q(\hat{\mathfrak{g}})$. Let $\pi$ be a $l$-tuple of polynomials. One can define the local Weyl module of the Yangian $Y(\mathfrak{g})$ similarly. We gave the definition of local Weyl module $W(\pi)$ by generators and defining relations in \cite{TaGu}, and proved that $W(\pi)$ is isomorphic to an ordered tensor product of fundamental representations when $\mathfrak{g}$ is classical.
The main challenge in \cite{TaGu} was to find an explicit cyclicity condition for an ordered tensor product, and the methodology will be introduced in Section 2 for more detail.

This work is a continuation of the paper \cite{TaGu}. We introduce a new algorithm to compute certain associated polynomials when $\mathfrak{g}$ is of type $G_2$, which enables us to provide an explicit cyclicity condition for an ordered tensor product of fundamental representations, see Theorem 1. Using this condition, we show that the local Weyl module $W(\pi)$ is isomorphic to an ordered tensor product of fundamental representations in Theorem 2.

\section{Preliminary}
In this section, we give a brief review of the previous paper \cite{TaGu} and one question left open in \textit{loc. cit.}

\begin{definition}
Let $\mathfrak{g}$ be a simple Lie algebra over $\mathbb{C}$ with rank $l$ and let $A=\left(a_{ij}\right)_{i,j\in I}$, where $I=\{1,2,\ldots, l\}$, be its Cartan matrix. Let $D=\operatorname{diag}\left(d_1,\ldots, d_{l}\right)$, $d_i\in \mathbb{N}$, such that $d_1, d_2,\ldots, d_l$ are co-prime and $DA$ is symmetric. The Yangian $Y(\frak g)$ is defined to be the
associative algebra with generators $x_{i,r}^{{}\pm{}}$,
$h_{i,r}$, $i\in I$, $r\in\mathbb{Z}_{\geq 0}$, and the following
defining relations:
\begin{equation*}\label{}
[h_{i,r},h_{j,s}]=0, \qquad [h_{i,0}, x_{j,s}^{\pm}]={}\pm
d_ia_{ij}x_{j,s}^{\pm}, \qquad [x_{i,r}^+, x_{j,s}^-]=\delta_{i,j}h_{i,r+s},
\end{equation*}
\begin{equation*}\label{}
[h_{i,r+1}, x_{j,s}^{\pm}]-[h_{i,r}, x_{j,s+1}^{\pm}]=
\pm\frac{1}{2}d_i
a_{ij}(h_{i,r}x_{j,s}^{\pm}+x_{j,s}^{\pm}h_{i,r}),
\end{equation*}
\begin{equation*}\label{}
[x_{i,r+1}^{\pm}, x_{j,s}^{\pm}]-
[x_{i,r}^{\pm}, x_{j,s+1}^{\pm}]=\pm\frac12
d_ia_{ij}(x_{i,r}^{\pm}x_{j,s}^{\pm}
+x_{j,s}^{\pm}x_{i,r}^{\pm}),
\end{equation*}
\begin{equation*}\label{}
\sum_\pi
[x_{i,r_{\pi(1)}}^{\pm},
[x_{i,r_{\pi(2)}}^{\pm}, \ldots,
[x_{i,r_{\pi(m)}}^{\pm},
x_{j,s}^{\pm}]\cdots]]=0, i\neq j,
\end{equation*}
for all sequences of non-negative integers $r_1,\ldots,r_m$, where
$m=1-a_{ij}$ and the sum is over all permutations $\pi$ of $\{1,\dots,m\}$.
\end{definition}

Denote by $V_m(a)$, $m\geq 1$ and $a\in \mathbb{C}$, the finite-dimensional irreducible representation of $Y(\mathfrak{sl}_2)$ associated to the Drinfeld polynomial $\big(u-a\big)\big(u-(a+1)\big)\ldots \big(u-(a+m-1)\big)$.

\begin{proposition}[Proposition 3.5,\cite{ChPr3}]\label{wra}
The module $V_m(a)$ has a basis
$\{w_0,w_1,\ldots, w_m\}$ on which the action of $Y(\mathfrak{sl}_2)$ is given by
\begin{center}
$ x_k^+w_s=(s+a)^k(s+1)w_{s+1},\ \ x_k^-w_s=(s+a-1)^k(m-s+1)w_{s-1},$\\

$h_kw_s=\big((s+a-1)^ks(m-s+1)-(s+a)^k(s+1)(m-s)\big)w_s$.
\end{center}
\end{proposition}

Suppose that a fixed reduced expression of the longest element of the Weyl group of $\mathfrak{g}$ is $w_0=s_{r_1}s_{r_2}\ldots s_{r_p}$, where $s_{r_j}$, for $1\leq j\leq p$, are simple reflections. Denote the $i$-th fundamental weight by $\omega_i$ and the $i$-th fundamental representation of $Y(\mathfrak{g})$ by $V_a(\omega_i)$. %, where $a\in \mathbb{C}$.
Let $v^{+}$ and $v^{-}$ be the highest and lowest weight vectors in the fundamental representation $V_a(\omega_i)$, respectively. Suppose $s_{r_{j+1}}s_{r_{j+2}}\ldots s_{r_p}(\omega_i)=m_j\omega_{r_j}+\sum\limits_{n\neq r_j}c_n\omega_n$. Then $$v^{-}=\big(x_{r_1,0}^{-}\big)^{m_{1}}\big(x_{r_2,0}^{-}\big)^{m_{2}}\ldots \big(x_{r_p,0}^{-}\big)^{m_{p}}v^{+}.$$
Define $s_{r_{j+1}}s_{r_{j+2}}\ldots s_{r_p}$ by $\sigma_j$ and $(x_{r_{j+1},0}^{-})^{m_{j+1}}(x_{r_{j+2},0}^{-})^{m_{j+2}}\ldots (x_{r_p,0}^{-})^{m_{p}}v^{+}$ by $v_{\sigma_j(\omega_i)}$. Let $Y_i$ be the subalgebra generated by $\{x_{i,r}^{\pm}, h_{i,r}|r\in\mathbb{Z}_{\geq 0}\}$ for $i\in I$, which is isomorphic to $Y(\mathfrak{sl}_2)$. Denote by $Y_{r_j}(v_{\sigma_j(\omega_{i})})$ the $Y_{r_j}$-module generated by the extremal vector $v_{\sigma_j(\omega_{i})}$.   We remark that it has been established in Section 5 in \cite{TaGu} that  $Y_{r_j}(v_{\sigma_j(\omega_{b_m})})$ is a highest weight representation and the degree of its associated polynomial is $m_j$. In \cite{TaGu}, we showed:

\noindent\textbf{Theorem 5.2} \cite{TaGu} An ordered tensor product $L=V_{a_1}(\omega_{b_1})\otimes V_{a_2}(\omega_{b_2})\otimes\ldots\otimes V_{a_k}(\omega_{b_k})$ is a highest weight representation if for all $1\leq j\leq p$ and $1\leq m<n\leq k$, when $b_n=r_j$, the difference of the number $\frac{a_n}{d_{r_j}}$ and any root of the associated polynomial of $Y_{r_j}(v_{\sigma_j(\omega_{b_m})})$ does not equal 1.

We proved the above theorem by adopting the ideas in \cite{Ch3}. To find an explicit cyclicity condition, it is enough to compute the associated polynomial of $Y_{r_j}(v_{\sigma_j(\omega_{b_m})})$ for $1\leq j\leq p$ and $1\leq m\leq l$.  When $\mathfrak{g}$ is a classical simple Lie algebra, by computing the associated polynomial of $Y_{r_j}(v_{\sigma_j(\omega_{b_1})})$ using some of the defining relations of $Y(\mathfrak{g})$, a concrete cyclicity condition of $L$ was obtained in Theorem 5.18 in \cite{TaGu}.

For the following reasons, it becomes complicated if one tries to describe a concrete cyclicity condition for the tensor product when $\mathfrak{g}$ is an exceptional simple Lie algebra . In this case, $m_j$ may be greater than or equal to 3. The computations of the eigenvalue of $h_{r_j, k}$ on the weight vector $v_{\sigma_j(\omega_i)}$, $3\leq k\leq m_{j}$, are tedious if one uses defining relations of $Y(\mathfrak{g})$.
In addition, the path from the highest weight vector $v^{+}$ to the lowest one $v^{-}$ is more subtle than in the case when $\mathfrak{g}$ is a classical simple Lie algebra. For instance, when $\mathfrak{g}$ is of type $F_4$, in $V_a(\omega_2)$, a path is:
\begin{eqnarray*}
% \nonumber to remove numbering (before each equation)
  v^{-} &=& x_{2,0}^{-}x_{1,0}^{-}(x_{3,0}^{-})^2(x_{2,0}^{-})^2x_{1,0}^{-}(x_{3,0}^{-})^2x_{2,0}^{-}(x_{4,0}^{-})^4 (x_{3,0}^{-})^4(x_{2,0}^{-})^{3}(x_{1,0}^{-})^3(x_{3,0}^{-})^2\\
   && \quad (x_{2,0}^{-})^2 (x_{3,0}^{-})^2(x_{4,0}^{-})^4(x_{3,0}^{-})^2 (x_{2,0}^{-})^2x_{1,0}^{-}(x_{3,0}^{-})^2x_{2,0}^{-}v^{+}.
\end{eqnarray*}
It seems to us that step-by-step computations are needed to compute the associated polynomial of $Y_{r_j}(v_{\sigma_j(\omega_{b_1})})$. The amount of time required to compute these associated polynomials would be tremendous using exactly the same approach as in \cite{TaGu}.

%The final goal is to shade some light to find the cyclicity condition for the tensor product of Kirillov-Reshetikin modules of $Y(\mathfrak{g})$, which will not be discussed in this paper.
\section{A concrete cyclicity condition and local Weyl modules}
From now on, let $\mathfrak{g}$ denote the simple Lie algebra of type $G_2$, unless the contrary is stated. Let $\alpha_1$ and $\alpha_2$ be the simple long and short roots, respectively (as labelled in \cite{Ca}) and let $\omega_1$ and $\omega_2$ be the fundamental weights.  The Cartan matrix of $\mathfrak{g}$ is $\left(
                          \begin{array}{cc}
                            2 & -1 \\
                            -3 & 2 \\
                          \end{array}
                        \right)$.
Let $D=\left(
       \begin{array}{cc}
         3 &\quad 0 \\
         0 &\quad 1 \\
       \end{array}
     \right),
$
so $DA$ is symmetric.
The Weyl group $\mathcal{W}$ of $\mathfrak{g}$ is generated by $s_1$ and $s_2$ such that $s_1(\alpha_1)=-\alpha_1,\ s_2(\alpha_1)=\alpha_1+3\alpha_2,\ s_1(\alpha_2)=\alpha_1+\alpha_2, \text{and}\ s_2(\alpha_2) = -\alpha_2.$
A reduced expression of the longest element in Weyl group is
$w_0=s_1s_2s_1s_2s_1s_2$. Using the method indicated in Section 2 of this paper, we found a path form $v^+$ to $v^{-}$.
\begin{lemma}Let $a\in \mathbb{C}$.
\begin{enumerate}
  \item In $V_a(\omega_1)$, $v^{-}=x_{1,0}^{-}(x_{2,0}^{-})^3(x_{1,0}^{-})^2(x_{2,0}^{-})^3x_{1,0}^{-}v^{+}.$
  \item In $V_a(\omega_2)$, $v^{-}=x_{2,0}^{-}x_{1,0}^{-}(x_{2,0}^{-})^2x_{1,0}^{-}x_{2,0}^{-}v^{+}$.
\end{enumerate}
\end{lemma}

Note that $x_{1,r}^{\pm}, h_{1,r}$ do not satisfy the defining relations of $Y(\mathfrak{sl}_2)$. Therefore, we need to re-scale the generators. Let $\tilde{x}_{1,r}^{\pm}=\frac{\sqrt{3}}{3^{r+1}}x_{1,r}^{\pm}$ and $\tilde{h}_{1,r}=\frac{1}{3^{r+1}}h_{1,r}$, then $\tilde{x}_{1,r}^{\pm}, \tilde{h}_{1,s}$ satisfy the defining relations of $Y(\mathfrak{sl}_2)$. Now we are in the position to compute the associated polynomial of $Y_{r_j}(v_{\sigma_j(\omega_{b_m})})$, where $b_m\in I=\{1,2\}$.
\begin{proposition}\label{g2fpv} Let $v^{+}\in V_a(\omega_1)$.\\
\begin{tabular}{|c|c|c|}
  \hline
  % after \\: \hline or \cline{col1-col2} \cline{col3-col4} ...
  Item & $Y(\mathfrak{sl}_2)$-module & Associated polynomial\\\hline
  1 & $Y_1(v^{+})$ & $u-\frac{a}{3}$ \\\hline
  2 & $Y_2(x_{1,0}^{-}v^{+})$ &$\big(u-(a+\frac{3}{2})\big)\big(u-(a+\frac{1}{2})\big)\big(u-(a-\frac{1}{2})\big)$ \\\hline
  3 & $Y_1\big((x_{2,0}^{-})^3x_{1,0}^{-}v^{+}\big)$ & $\big(u-\frac{a+2}{3}\big)\big(u-\frac{a+1}{3}\big)$ \\\hline
  4 & $Y_2\big((x_{1,0}^{-})^2(x_{2,0}^{-})^3x_{1,0}^{-}v^{+}\big)$ & $\big(u-(a+\frac{7}{2})\big)\big(u-(a+\frac{5}{2})\big)\big(u-(a+\frac{3}{2})\big)$\\\hline
  5 & $Y_1\big((x_{2,0}^{-})^3(x_{1,0}^{-})^2(x_{2,0}^{-})^3x_{1,0}^{-}v^{+}\big)$ & $u-\frac{a+3}{3}$ \\\hline
\end{tabular}
\end{proposition}
\begin{proof}
We omit the proof of items  1, 2 and 5 since these proofs can be checked using the same approach as in paper \cite{TaGu}. Lemma \ref{c7c1l2} is devoted to proving the third item. The fourth item is proved in Lemma \ref{c7c1l3}.
\end{proof}
Let $H$ be the subalgebra of $Y(\mathfrak{g})$ generated by all $h_{i,k}$ and $h_i(u)=1+h_{i,0}u^{-1}+h_{i,1}u^{-2}+\ldots$, where $i\in I=\{1,2\}$ and $k\in \mathbb{Z}_{\geq 0}$. We use alternate generators for $H$ as given in \cite{Le}. Let $$H_{i}(u)=\sum\limits_{k=0}^{\infty} H_{i,k}u^{-k-1}:=\ln\big(h_i(u)\big).$$
An explicit computation shows that
\begin{align}\label{equ21}
H_{i}(u)=h_{i,0}u^{-1}+\big(h_{i,1}-\frac{1}{2}(h_{i,0})^2\big)u^{-2}
+\big(h_{i,2}-h_{i,0}h_{i,1}+\frac{1}{3}(h_{i,0})^3\big)u^{-3}+\ldots.
\end{align}

\begin{lemma}[Corollary 1.5, \cite{Le}]\label{c7lecor15} Let $\mathfrak{g}$ be a complex simple Lie algebra.
\begin{align*}
[H_{i,k}, x_{j,l}^{\pm}]=&\pm d_ia_{ij}x_{j,l+k}^{\pm}\pm\sum_{\substack{0\leq s\leq k-2\\ k+s\ \text{even}}}2^{s-k}(d_ia_{ij})^{k+1-s}\frac{{ k+1\choose s}}{k+1}x_{j,l+s}^{\pm}.
\end{align*}
\end{lemma}
%\begin{lemma}[Corollary 1.5, \cite{Le}]\label{c7lecor15} When $i\neq j$,
%\begin{align*}
%[H_{i,k}, x_{j,l}^{\pm}]=&\mp 3x_{j,l+k}^{\pm}\pm \sum_{\substack{0\leq s\leq k-2\\ k+s\ \text{even}}}2^{s-k}(-3)^{k+1-s}\frac{{ k+1\choose s}}{k+1}x_{j,l+s}^{\pm}.
%\end{align*}
%\end{lemma}

Let $V$ be a finite-dimensional highest weight representation of $Y(\mathfrak{sl}_2)$ whose associated polynomial is $\pi$. Let $v^{+}$ and $v^{-}$ be highest and lowest weight vectors of $V$, respectively.
%We denote the generators of $Y(\mathfrak{sl}_2)$ by $x_{k}^{\pm}$ and $h_{k}$.
\begin{lemma}\label{huv-iqg2}
$h(u)v^{-}=\frac{\pi(u-1)}{\pi(u)}v^{-}$.
\end{lemma}
\begin{proof}
It is enough to consider the case when $V$ is irreducible. It was established in Proposition 3.1 \cite{ChPr4} that $h(u)v^{-}=\frac{p(u)}{p(u+1)}v^{-}$, where $p(u)$ is the associated polynomial of the right dual $V^t$. It was showed  in Proposition 2.4 in \cite{ChPr3} that $p(u)=\pi(u-1)$. Therefore this lemma is proved. 
\end{proof}
To show Lemma \ref{c7c1l2}, we need the following corollary. We remark that in the proof of the corollary, we use some algorithm which did not show up in \cite{TaGu}.
\begin{corollary}\label{g2c21v} In the representation $Y_2(x_{1,0}^{-}v^{+})$,
\begin{align*}
1.\qquad \big(x_{2,1}^{-}(x_{2,0}^{-})^2&+x_{2,0}^{-}x_{2,1}^{-}x_{2,0}^{-}+(x_{2,0}^{-})^2 x_{2,1}^{-}\big)x_{1,0}^{-}v^{+}\\
       &=\Big((a+\frac{3}{2})+(a+\frac{1}{2})+(a-\frac{1}{2})\Big)(x_{2,0}^{-})^3x_{1,0}^{-}v^{+}.\\
2.\qquad \big(x_{2,2}^{-}(x_{2,0}^{-})^2&+x_{2,0}^{-}x_{2,2}^{-}x_{2,0}^{-}+(x_{2,0}^{-})^2 x_{2,2}^{-}\big)x_{1,0}^{-}v^{+}\\
      &=\Big((a+\frac{3}{2})^2+(a+\frac{1}{2})^2+
      (a-\frac{1}{2})^2\Big)(x_{2,0}^{-})^3x_{1,0}^{-}v^{+}.
\end{align*}
\end{corollary}
\begin{proof}
It follows from Proposition \ref{g2fpv} that the associated polynomial of the highest weight representation $Y_2(x_{1,0}^{-}v^{+})$ is $\big(u-(a+\frac{3}{2})\big)\big(u-(a+\frac{1}{2})\big)\big(u-(a-\frac{1}{2})\big)$. Therefore, $h_{2}(u)(x_{1,0}^{-}v^{+})=\frac{u-(a-\frac{3}{2})}{u-(a+\frac{3}{2})}
x_{1,0}^{-}v^{+}$ and $h_{2}(u)\big((x_{2,0}^{-})^3x_{1,0}^{-}v^{+}\big)=\frac{u-(a+\frac{5}{2})}{u-(a-\frac{1}{2})}
(x_{2,0}^{-})^3x_{1,0}^{-}v^{+}$ by Lemma \ref{huv-iqg2}.
Thus $$H_{2}(u)(x_{1,0}^{-}v^{+})=\Big(\ln(1-(a-\frac{3}{2})u^{-1})-\ln(1-(a+\frac{3}{2})u^{-1})\Big)
x_{1,0}^{-}v^{+}$$ and $$H_{2}(u)\big((x_{2,0}^{-})^3x_{1,0}^{-}v^{+}\big)=\Big(\ln(1-(a+\frac{5}{2})u^{-1})-\ln(1-(a-\frac{1}{2})u^{-1})\Big)
(x_{2,0}^{-})^3x_{1,0}^{-}v^{+}.$$
In particular, we have both $H_{2,1}(x_{2,0}^{-})^3x_{1,0}^{-}v^{+}=\frac{1}{2}\big((a-\frac{1}{2})^2-(a+\frac{5}{2})^2\big)
(x_{2,0}^{-})^3x_{1,0}^{-}v^{+}$ and $H_{2,2}(x_{2,0}^{-})^3x_{1,0}^{-}v^{+}=\frac{1}{3}\big((a-\frac{1}{2})^3-(a+\frac{5}{2})^3\big)
(x_{2,0}^{-})^3x_{1,0}^{-}v^{+}$. By Lemma \ref{c7lecor15}, $[H_{2,1}, x_{2,0}^{-}]=-2x_{2,1}^{-}$ and $[H_{2,2}, x_{2,0}^{-}]=-2x_{2,2}^{-}-\frac{2}{3}x_{2,0}^{-}$. We are going to show the first item in this corollary.
\begin{align*}
H_{2,1}&(x_{2,0}^{-})^3x_{1,0}^{-}v^{+}\\
&=[H_{2,1},(x_{2,0}^{-})^3]x_{1,0}^{-}v^{+}+(x_{2,0}^{-})^3 H_{2,1}x_{1,0}^{-}v^{+}\\
&=-2\big(x_{2,1}^{-}(x_{2,0}^{-})^2+x_{2,0}^{-}x_{2,1}^{-}x_{2,0}^{-}+(x_{2,0}^{-})^2 x_{2,1}^{-}\big)x_{1,0}^{-}v^{+}\\
&+\frac{1}{2}\big((a+\frac{3}{2})^2-(a-\frac{3}{2})^2\big)(x_{2,0}^{-})^3 x_{1,0}^{-}v^{+}.
\end{align*}
Therefore,
\begin{align*}
\big(x_{2,1}^{-}(x_{2,0}^{-})^2&+x_{2,0}^{-}x_{2,1}^{-}x_{2,0}^{-}+(x_{2,0}^{-})^2 x_{2,1}^{-}\big)x_{1,0}^{-}v^{+}\\
&=\frac{1}{4}\big((a+\frac{5}{2})^2-(a-\frac{1}{2})^2+(a+\frac{3}{2})^2-(a-\frac{3}{2})^2\big)(x_{2,0}^{-})^3 x_{1,0}^{-}v^{+}\\
&=(3a+\frac{3}{2})(x_{2,0}^{-})^3x_{1,0}^{-}v^{+}.
\end{align*}
The second item in this corollary can be obtained similarly, so we omit the proof. 
\end{proof}
\begin{lemma}\label{c7c1l2}
The associated polynomial of the representation $Y_1\Big((x_{2,0}^{-})^3x_{1,0}^{-}v^{+}\Big)$ is given by $\big(u-\frac{a+2}{3}\big)\big(u-\frac{a+1}{3}\big)$.
\end{lemma}
\begin{proof} The associated polynomial of the representation $Y_1\Big((x_{2,0}^{-})^3x_{1,0}^{-}v^{+}\Big)$ is of degree 2, say $(u-a_1)(u-a_2)$. The eigenvalues of $(x_{2,0}^{-})^3x_{1,0}^{-}v^{+}$ under $\tilde{h}_{1,1}$ and $\tilde{h}_{1,2}$ will tell us the values of $a_1$ and $a_2$. We first compute the eigenvalues of $(x_{2,0}^{-})^3x_{1,0}^{-}v^{+}$ under $H_{1,1}$ and $H_{1,2}$.
\begin{align*}
H_{1,1}&(x_{2,0}^{-})^3x_{1,0}^{-}v^{+}\\
&=[H_{1,1},x_{2,0}^{-}](x_{2,0}^{-})^2x_{1,0}^{-}v^{+}+x_{2,0}^{-}[H_{1,1},x_{2,0}^{-}]x_{2,0}^{-}x_{1,0}^{-}v^{+}\\
&\qquad +(x_{2,0}^{-})^2[H_{1,1},x_{2,0}^{-}]x_{1,0}^{-}v^{+}+(x_{2,0}^{-})^3 H_{1,1}x_{1,0}^{-}v^{+}\\
&=3x_{2,1}^{-}(x_{2,0}^{-})^2x_{1,0}^{-}v^{+}+3x_{2,0}^{-}x_{2,1}^{-}x_{2,0}^{-}x_{1,0}^{-}v^{+}+3(x_{2,0}^{-})^2x_{2,1}^{-}x_{1,0}^{-}v^{+}\\
&\qquad+\frac{1}{2}\big((a-3)^2-a^2\big)(x_{2,0}^{-})^3 x_{1,0}^{-}v^{+}\\
&=3\big((a-\frac{1}{2})+(a+\frac{1}{2})+(a+\frac{3}{2})-(a+\frac{3}{2})\big)
(x_{2,0}^{-})^3x_{1,0}^{-}v^{+}\\
&=6a(x_{2,0}^{-})^3x_{1,0}^{-}v^{+},
\end{align*}where the third equality follows from the first item of Corollary \ref{g2c21v}.

\begin{align*}
&H_{1,2}(x_{2,0}^{-})^3x_{1,0}^{-}v^{+}\\
&=[H_{1,2},x_{2,0}^{-}](x_{2,0}^{-})^2x_{1,0}^{-}v^{+}+x_{2,0}^{-}[H_{1,2},x_{2,0}^{-}]x_{2,0}^{-}x_{1,0}^{-}v^{+}\\
&\qquad +(x_{2,0}^{-})^2[H_{1,2},x_{2,0}^{-}]x_{1,0}^{-}v^{+}+(x_{2,0}^{-})^3 H_{1,2}x_{1,0}^{-}v^{+}\\
&=\big(3x_{2,2}^{-}+(\frac{3}{2})^2x_{2,0}^{-}\big)(x_{2,0}^{-})^2x_{1,0}^{-}v^{+}+
x_{2,0}^{-}\big(3x_{2,2}^{-}+(\frac{3}{2})^2x_{2,0}^{-}\big)x_{2,0}^{-}x_{1,0}^{-}v^{+}\\
&+(x_{2,0}^{-})^2\big(3x_{2,2}^{-}+(\frac{3}{2})^2x_{2,0}^{-}\big)x_{1,0}^{-}v^{+}
-\frac{1}{3}\big((a-3)^3-a^3\big)(x_{2,0}^{-})^3 x_{1,0}^{-}v^{+}\\
&=3\big((a-\frac{1}{2})^2+(a+\frac{1}{2})^2+(a+\frac{3}{2})^2+\frac{9}{4}-(a^2+3a+3)\big)(x_{2,0}^{-})^3x_{1,0}^{-}v^{+}\\
&=(6a^2+6)(x_{2,0}^{-})^3x_{1,0}^{-}v^{+}.
\end{align*}
It follows from equation (\ref{equ21}) and the above computations that
\begin{center}
$\tilde{h}_{1,1}(x_{2,0}^{-})^3x_{1,0}^{-}v^{+}=(\frac{2a}{3}+2)(x_{2,0}^{-})^3x_{1,0}^{-}v^{+};$
$\tilde{h}_{1,2}(x_{2,0}^{-})^3x_{1,0}^{-}v^{+}=(2(\frac{a}{3})^2+\frac{4 a}{3}+\frac{14}{9})(x_{2,0}^{-})^3x_{1,0}^{-}v^{+}.
$
\end{center}

It follows from relation (5.1) in \cite{TaGu} that $a_1+a_2=\frac{2a}{3}+1$ and $a_1^2+a_2^2+a_1+a_2=2(\frac{a}{3})^2+\frac{4 a}{3}+\frac{14}{9}$. Then $a_1=\frac{a+1}{3}$ and $a_2=\frac{a+2}{3}$, or vice-versa with $a_1$ and $a_2$ switched. Thus the associated polynomial of the representation $Y_1\Big((x_{2,0}^{-})^3x_{1,0}^{-}v^{+}\Big)$ is $\big(u-\frac{a+2}{3}\big)\big(u-\frac{a+1}{3}\big)$.
\end{proof}
Similar to the proof of Corollary \ref{g2c21v}, we can prove the following corollary. We remark that:
$h_1(u)\big((x_{2,0}^{-})^3x_{1,0}^{-}v^{+}\big)=\frac{u-(a+1)}{u-(a+2)}\frac{u-a}{u-(a+1)}(x_{2,0}^{-})^3
    x_{1,0}^{-}v^{+}$.
\begin{corollary}\label{c7c1r2} In the representation $Y_1\Big((x_{2,0}^{-})^3x_{1,0}^{-}v^{+}\Big)$,
\begin{enumerate}
\item $(x_{1,1}^{-}x_{1,0}^{-}+x_{1,0}^{-}x_{1,1}^{-})(x_{2,0}^{-})^3
    x_{1,0}^{-}v^{+}=\big((a+1)+(a+2)\big)(x_{1,0}^{-})^2(x_{2,0}^{-})^3x_{1,0}^{-}v^{+}$;
\item $(x_{1,2}^{-}x_{1,0}^{-}+x_{1,0}^{-}x_{1,2}^{-})(x_{2,0}^{-})^3
    x_{1,0}^{-}v^{+}=\big((a+1)^2+(a+2)^2\big)(x_{1,0}^{-})^2(x_{2,0}^{-})^3x_{1,0}^{-}v^{+}$;
\item $(x_{1,3}^{-}x_{1,0}^{-}+x_{1,0}^{-}x_{1,3}^{-})(x_{2,0}^{-})^3x_{1,0}^{-}v^{+}=\big((a+1)^3+(a+2)^3\big)(x_{1,0}^{-})^2(x_{2,0}^{-})^3x_{1,0}^{-}v^{+}$.
\end{enumerate}
\end{corollary}

\begin{lemma}\label{c7c1l3}
The associated polynomial of the module $Y_2\big((x_{1,0}^{-})^2(x_{2,0}^{-})^3x_{1,0}^{-}v^{+}\big)$ is given by $\big(u-(a+\frac{3}{2})\big)\big(u-(a+\frac{5}{2})\big)\big(u-(a+\frac{7}{2})\big).$
\end{lemma}
\begin{proof} The associated polynomial of the representation $Y_2\big((x_{1,0}^{-})^2(x_{2,0}^{-})^3x_{1,0}^{-}v^{+}\big)$ is of degree 3, say $(u-a_1)(u-a_2)(u-a_3)$. The eigenvalues of $(x_{1,0}^{-})^2(x_{2,0}^{-})^3x_{1,0}^{-}v^{+}$ under $h_{2,1}$, $h_{2,2}$ and $h_{2,3}$ will tell us the values of $a_1$, $a_2$ and $a_3$.
\begin{align*}
% \nonumber to remove numbering (before each equation)
H_{2,1}&(x_{1,0}^{-})^2(x_{2,0}^{-})^3x_{1,0}^{-}v^{+} \\
&=[H_{2,1},x_{1,0}^{-}]x_{1,0}^{-}(x_{2,0}^{-})^3x_{1,0}^{-}v^{+}+
x_{1,0}^{-}[H_{2,1},x_{1,0}^{-}](x_{2,0}^{-})^3x_{1,0}^{-}v^{+}+(x_{1,0}^{-})^2 H_{2,1}(x_{2,0}^{-})^3x_{1,0}^{-}v^{+}\\
&=3\big(x_{1,1}^{-}x_{1,0}^{-}+x_{1,0}^{-}x_{1,1}^{-}\big)(x_{2,0}^{-})^3x_{1,0}^{-}v^{+}
 -\frac{1}{2}\big((a+\frac{5}{2})^2-(a-\frac{1}{2})^2\big)(x_{1,0}^{-})^2(x_{2,0}^{-})^3x_{1,0}^{-}v^{+}\\
&=\big((6a+9)-(3a+3)\big)(x_{1,0}^{-})^2(x_{2,0}^{-})^3x_{1,0}^{-}v^{+}\\
&=(3a+6)(x_{1,0}^{-})^2(x_{2,0}^{-})^3x_{1,0}^{-}v^{+},
\end{align*}
where the third equality follows from the first item of Corollary \ref{c7c1r2}.
\begin{align*}
% \nonumber to remove numbering (before each equation)
H_{2,2}& (x_{1,0}^{-})^2(x_{2,0}^{-})^3x_{1,0}^{-}v^{+} \\
&=[H_{2,2},x_{1,0}^{-}]x_{1,0}^{-}(x_{2,0}^{-})^3x_{1,0}^{-}v^{+}+
x_{1,0}^{-}[H_{2,2},x_{1,0}^{-}](x_{2,0}^{-})^3x_{1,0}^{-}v^{+}\\
&\qquad +(x_{1,0}^{-})^2 H_{2,2}(x_{2,0}^{-})^3x_{1,0}^{-}v^{+}\\
&=(3x_{1,2}^{-}+\frac{9}{4}x_{1,0}^{-})x_{1,0}^{-}(x_{2,0}^{-})^3x_{1,0}^{-}v^{+}+x_{1,0}^{-}(3x_{1,2}^{-}+\frac{9}{4}x_{1,0}^{-})(x_{2,0}^{-})^3x_{1,0}^{-}v^{+}\\
&\qquad-\frac{1}{3}\big((a+\frac{5}{2})^3-(a-\frac{1}{2})^3\big)(x_{1,0}^{-})^2
(x_{2,0}^{-})^3x_{1,0}^{-}v^{+}\\
&=(3a^2+12a+\frac{57}{4})(x_{1,0}^{-})^2(x_{2,0}^{-})^3x_{1,0}^{-}v^{+}.
\end{align*}
\begin{align*}
% \nonumber to remove numbering (before each equation)
H_{2,3}&(x_{1,0}^{-})^2(x_{2,0}^{-})^3x_{1,0}^{-}v^{+} \\
&=[H_{2,3},x_{1,0}^{-}]x_{1,0}^{-}(x_{2,0}^{-})^3x_{1,0}^{-}v^{+}+
x_{1,0}^{-}[H_{2,3},x_{1,0}^{-}](x_{2,0}^{-})^3x_{1,0}^{-}v^{+}\\
&\qquad +(x_{1,0}^{-})^2 H_{2,3}(x_{2,0}^{-})^3x_{1,0}^{-}v^{+}\\
&=(3x_{1,3}^{-}+\frac{27}{4}x_{1,1}^{-})x_{1,0}^{-}(x_{2,0}^{-})^3x_{1,0}^{-}v^{+}
+x_{1,0}^{-}(3x_{1,3}^{-}+\frac{27}{4}x_{1,1}^{-})(x_{2,0}^{-})^3x_{1,0}^{-}v^{+}\\
&\qquad-\frac{1}{4}\big((a+\frac{5}{2})^4-(a-\frac{1}{2})^4\big)(x_{1,0}^{-})^2(x_{2,0}^{-})^3x_{1,0}^{-}v^{+}\\
&=(3a^3+18a^2+\frac{171}{4}a+\frac{75}{2})(x_{1,0}^{-})^2(x_{2,0}^{-})^3x_{1,0}^{-}v^{+}.
\end{align*}
It follows from equation (\ref{equ21}) and the above computations that
\begin{align*}
h_{2,1}(x_{1,0}^{-})^2(x_{2,0}^{-})^3x_{1,0}^{-}v^{+}&=
3(a+\frac{7}{2})(x_{1,0}^{-})^2(x_{2,0}^{-})^3x_{1,0}^{-}v^{+};\\
h_{2,2}(x_{1,0}^{-})^2(x_{2,0}^{-})^3x_{1,0}^{-}v^{+}&=
3(a+\frac{7}{2})^2(x_{1,0}^{-})^2(x_{2,0}^{-})^3x_{1,0}^{-}v^{+};\\
h_{2,3}(x_{1,0}^{-})^2(x_{2,0}^{-})^3x_{1,0}^{-}v^{+}&=
3(a+\frac{7}{2})^3(x_{1,0}^{-})^2(x_{2,0}^{-})^3x_{1,0}^{-}v^{+}.
\end{align*}
Note that \begin{align*}
&\frac{u-(a_1-1)}{u-a_1}\frac{u-(a_2-1)}{u-a_2} \frac{u-(a_3-1)}{u-a_3}\nonumber\\
&\quad=1+3u^{-1}+(a_1+a_2+a_3+3)u^{-2}+(a_1^2+a_2^2+a_3^2+2a_1+2a_2+2a_3+1)u^{-3}\\
&\qquad +(a_1^3+a_2^3+a_3^3+2a_1^2+2a_2^2+2a_3^2+a_1a_2+a_1a_3+a_2a_3+a_1+a_2+a_3)u^{-4}\\
&\qquad+\ldots.
\end{align*}
Computations show that $a_1=a+\frac{3}{2}$, $a_2=a+\frac{5}{2}$ and $a_3=a+\frac{7}{2}$, or vice-versa with $a_1$, $a_2$ and $a_3$ switched. Thus the associated polynomial of the representation $Y_2((x_{1,0}^{-})^2(x_{2,0}^{-})^3x_{1,0}^{-}v^{+})$ is given by $\big(u-(a+\frac{3}{2})\big)\big(u-(a+\frac{5}{2})\big)\big(u-(a+\frac{7}{2})\big).$ 
\end{proof}
We will omit the proof of the following proposition since it can be proved using some of the defining relations of $Y(\mathfrak{g})$.
\begin{proposition}Let $v^{+}\in V_a(\omega_2)$.\\
\begin{tabular}{|c|c|c|}
  \hline
  % after \\: \hline or \cline{col1-col2} \cline{col3-col4} ...
  Item & $Y(\mathfrak{sl}_2)$-module & Associated polynomial\\\hline
  1 & $Y_2(v^{+})$ & $u-a$ \\\hline
  2 & $Y_1(x_{2,0}^{-}v^{+})$ & $u-(\frac{a}{3}+\frac{1}{2})$ \\\hline
  3 & $Y_2\Big(x_{1,0}^{-}x_{2,0}^{-}v^{+}\Big)$ & $\big(u-(a+3)\big)\big(u-(a+2)\big)$ \\\hline
  4 & $Y_1((x_{2,0}^{-})^2x_{1,0}^{-}x_{2,0}^{-}v^{+})$ & $u-(\frac{a}{3}+\frac{7}{6})$ \\\hline
  5 & $Y_2(x_{1,0}^{-}(x_{2,0}^{-})^2x_{1,0}^{-}x_{2,0}^{-}v^{+})$ & $u-(a+5)$ \\\hline
\end{tabular}
\end{proposition}
We summarize all results in Propositions 1 and 2 into the coming corollary. Denote by $T(b_m, r_j)$ the set of all possible roots of the associated polynomial of $Y_{r_j}(v_{\sigma_j(\omega_{b_m})})$.
\begin{corollary}
 $T(1,1)=\{\frac{a}{3},\frac{a+1}{3},\frac{a+2}{3},\frac{a+3}{3}\}$,
$T(1,2)=\{a-\frac{1}{2}, a+\frac{1}{2}, a+\frac{3}{2},a+\frac{5}{2},a+\frac{7}{2}\}$, $T(2,1)=\{\frac{a}{3}+\frac{1}{2}, \frac{a}{3}+\frac{7}{6}\}$ and $T(2,2)=\{a, a+2, a+3, a+5\}$.
\end{corollary}
By Theorem 5.2 in \cite{TaGu} and Proposition 3.8 in \cite{ChPr8}, we have the following theorem.
\begin{theorem}\label{g2clihwr}
Let $L=V_{a_1}(\omega_{b_1})\otimes V_{a_2}(\omega_{b_2})\otimes\ldots\otimes V_{a_k}(\omega_{b_k})$ be an ordered tensor product of fundamental representations of $Y(\mathfrak{g})$, and define $S(1,1)=\{3,4,5,6\}$, $S(1,2)=\{\frac{1}{2},\frac{3}{2},\frac{5}{2},\frac{7}{2},\frac{9}{2}\}$, $S(2,1)=\{\frac{9}{2},\frac{13}{2}\}$ and $S(2,2)=\{1,3,4,6\}$.
\begin{enumerate}
  \item If $a_j-a_i\notin S(b_i, b_j)$ for $1\leq i<j\leq k$, then $L$ is a highest weight representation of $Y(\mathfrak{g})$.
  \item If $a_j-a_i\notin S(b_i, b_j)$ for $1\leq i\neq j\leq k$, then $L$ is an irreducible representation of $Y(\mathfrak{g})$.
\end{enumerate}
\end{theorem}
\begin{remark}\

\begin{itemize}
  \item[(i).] In Section 6.2 \cite{Ch3}, V. Chari gave the set $\mathrm{S}(i_1, i_2)$, $i_1\leq i_2$ of values of $a_1^{-1}a_2$ for which the tensor product $V(i_1, a_1)\otimes V(i_2, a_2)$ may fail to be irreducible as a module over the quantum loop algebra. Note that, in the following discussion, we interchange the node labels on the Dynkin diagram of $G_2$ as used in that paper.
V. Chari found that
\begin{align*}
\mathrm{S}(1,1)=\{q^6, q^8, q^{10}, q^{12}\}.\qquad\qquad
&\mathrm{S}(2,1)=\{q^7, q^{11}\}.\\
\mathrm{S}(1,2)=\{q^3, q^7\}.\qquad\qquad\qquad\quad\
&\mathrm{S}(2,2)=\{q^2, q^6, q^{8}, q^{12}\}.
\end{align*}

  \item[(ii).] In this paragraph, let $\mathfrak{g}$ be any finite-dimensional complex simple Lie algebra of rank $l$. The fundamental representations of $Y(\mathfrak{g})$ can be treated as special cases of the Kirillov-Reshetikhin modules, which are the finite-dimensional irreducible representations associated to an $l$-tuple of polynomials $\operatorname{\pi}=\big(\pi_1(u), \pi_2(u),\ldots, \pi_l(u)\big)$ such that $\pi_i(u)=\big(1-au^{-1}\big)\big(1-(a+1)u^{-1}\big)\ldots \big(1-(a+m-1)u^{-1}\big)$ and $\pi_j(u)=1$, for all $j\neq i$. The methods used in this paper could shed some light on obtaining a concrete cyclicity condition for the tensor product of Kirillov-Reshetikhin modules of $Y(\mathfrak{g})$.
\end{itemize}
\end{remark}

We close this section by providing the structure of $W(\pi)$. To obtain an upper bound on the dimension of $W(\pi)$, we use the dimension of the local Weyl module $W(\lambda)$ of the current algebra $\mathfrak{g}[t]$, which is given in \cite{Na}.
\begin{proposition}[Corollary 9.5, \cite{Na}]\label{dwmocsp}
Let $\lambda= m_1\omega_1+m_2\omega_2$. Then
$$\operatorname{Dim}\big(W(\lambda)\big)=\Big(\operatorname{Dim}\big(W(\omega_1)\big)\Big)^{m_1}\Big(\operatorname{Dim}\big(W(\omega_2)\big)\Big)^{m_2}.$$
\end{proposition}

It follows from Theorem \ref{g2clihwr} and Proposition 2.15 in \cite{ChPr4} that
\begin{proposition}\label{c7mtotsolwm} Let $\pi=\big(\pi_1(u), \pi_{2}(u)\big)$ be a pair of monic polynomials in $u$, and let $\pi_i(u)=\prod\limits_{j=1}^{m_i} (u-a_{i,j})$.  Let $S=\{a_{1,1},\ldots, a_{1,m_1}, a_{2,1}, \ldots, a_{2,m_2}\}$ be a multi-set of roots. Let  $a_1=a_{i,j}$ be one of the numbers in $S$ with maximal real part and let $b_1=i$. Similarly, let $a_r=a_{s,t}(r\geq 2)$ be one of the numbers in $S\setminus\{a_1, \ldots, a_{r-1}\}$ with maximal real part and let $b_r=s$.  Let $L=V_{a_1}(\omega_{b_1})\otimes V_{a_2}(\omega_{b_2})\otimes\ldots\otimes V_{a_k}(\omega_{b_k})$, where $k=m_1+m_2$. Then $L$ is a highest weight representation and its associated polynomials are $\pi_1(u)$ and $\pi_{2}(u)$. %$A=\{a\in \mathbb{C}| \text{there exists $p_i$ such that $p_i(a)=0$}\}$
\end{proposition}

\begin{theorem}The local Weyl module $W(\pi)$ of $Y(\mathfrak{g})$ associated to $\pi$ is isomorphic to $L$ as in Proposition \ref{c7mtotsolwm}.
\end{theorem}
\begin{proof}
Let $\lambda= m_1\omega_1+m_2\omega_2$. On the one hand, $\operatorname{Dim}\big(W(\pi)\big)\leq \operatorname{Dim}\big(W(\lambda)\big)$ by Theorem 3.8 \cite{TaGu}; on the other
hand, $\operatorname{Dim}\big(W(\pi)\big)\geq \operatorname{Dim}(L)$ by the maximality of the local Weyl modules of Yangians. Note that as $G_2$-modules, $W(\omega_i)\cong KR(\omega_i)\cong V_a(\omega_i)$ and the latter isomorphism follows easily from the main theorem of Section 2.3 in \cite{ChMo} and Theorem 6.3 in \cite{ChPr4}. In particular, $\operatorname{Dim}\big(W(\omega_i)\big)=\operatorname{Dim}\big(V_{a}(\omega_{i})\big)$ for any $a\in \mathbb{C}$. Therefore,
$\operatorname{Dim}\big(W(\lambda)\big)=\operatorname{Dim}(L)$. This implies that
$\operatorname{Dim}\big(W(\pi)\big)=\operatorname{Dim}(L)$, and therefore $W(\pi)\cong L$. 
\end{proof}
\section*{Acknowledgement}
The author's research was supported by the NSERC Discovery grant of his supervisor, N. Guay. The author would like to thank the referee for the valuable comments and suggestions.

\end{document}